\documentclass[10pt,twoside]{amsart}
  \usepackage{amsmath,amssymb,amsthm,amscd,latexsym,graphicx,tikz}
 \usepackage{epigraph}
 \usepackage[OT2,OT1]{fontenc}
  \usepackage{enumerate}

\usepackage{graphicx}
\usepackage[all]{xy}

\usepackage{graphicx}
\usepackage[all]{xy}

 \makeatletter
 \@specialpagefalse\headheight=8.5pt
 \makeatother

 \renewcommand{\to}{\longrightarrow}

 \newcommand{\A}{\mathbb A}

 \newcommand{\F}{\mathbb{F}}
 \newcommand{\G}{\mathbb{G}}
 \newcommand{\Bl}{\mathrm{Bl}}
 \newcommand{\Gm}{\mathbb{G}_{\mathrm {m}}}

 \newcommand{\N}{\mathbb{N}}
 \newcommand{\Pic}{\mathrm{Pic}}
 \renewcommand{\P}{\mathbb P}
 \newcommand{\Q}{\mathbb{Q}}
 \newcommand{\R}{\mathbb{R}}
 
 \newcommand{\Z}{\mathbb{Z}}

 \newcommand{\CH}{\mathrm{CH} } 
 \newcommand{\Stab}{{\mathrm{Stab}}}

 \newcommand{\GL}{\mathrm{GL}}

  \newcommand{\SStab}{\mathbf{Stab}}

 \newcommand{\AAut}{\mathbf{Aut}}
 \newcommand{\Aut}{\mathrm{Aut}}

 \newcommand{\carac}{\mathrm{char}}

  \newcommand{\End}{\mathrm{End}}
 
 \newcommand{\Proj}{\mathrm{Proj}}

 \newcommand{\Hom}{\mathrm{Hom}}

  \newcommand{\Ver}{\mathrm{Ver}}

 \newcommand{\Lie}{\mathrm{Lie}}

\newcommand{\Id}{\mathrm{Id}}

\renewcommand{\Im}{\mathrm{Im}}

 \newcommand{\PGL}{\mathrm{PGL}}
 
\newcommand{\Ker}{\mathrm{Ker}}
 
 \newcommand{\Spec}{\mathrm{Spec}}

 \newcommand{\Gr}{\mathrm{Gr}}

\newcommand{\Sym}{\mathrm{Sym}}

 \theoremstyle{plain}
 \newtheorem{thm}{Theorem}[section]
 \newtheorem{defi}[thm]{Definition}

 \newtheorem{prop}[thm]{Proposition}

 \newtheorem{lem}[thm]{Lemma}
 \newtheorem{coro}[thm]{Corollary}

 \theoremstyle{remark}
 \newtheorem{rem}[thm]{Remark}
 \newtheorem{qu}[thm]{Question}

 \newtheorem{ex}[thm]{Example}

 \newenvironment{dem}{{\bf Proof.}}{\hfill$\square$}

\date{\today}
 
\begin{document}

\title{ Realisation of linear algebraic groups as automorphism groups}

\author{Mathieu Florence}
\address{Sorbonne Université and Université Paris Cité, CNRS, IMJ-PRG, F-75005 Paris, France.}
\email{mathieu.florence@imj-prg.fr}




\begin{abstract}

Let $G$ be a linear algebraic group, over a field $F$. We show that $G$ is isomorphic to the automorphism group scheme of a smooth projective $F$-variety, defined as the blow-up of a projective space, along a suitable smooth subvariety.
\end{abstract}
\maketitle
\newpage
\tableofcontents
\newpage
\section{Introduction.}
Let  $X$ be a projective variety over a field $F$. The automorphism group functor $\AAut(X)$ is represented by a group scheme, locally of finite type over $F$. This  is due to Grothendieck (see also \cite{MO}, Theorem 3.7). Note that the sub-group scheme $\AAut^0(X) \subset \AAut(X)$, defined as the connected component of the identity, is then a group scheme of finite type over $F$- that is to say, an algebraic group over $F$.\\ Conversely, it is natural to ask:

\begin{qu}
 Let  $G$ be an algebraic group over a field $F$.\\
Does there exists a smooth projective $F$-variety $X$, such that $G \simeq \AAut(X)$? 
\end{qu}

When $G=A$ is an abelian variety, the answer was found independently by several authors:  it is positive, if and only if $\Aut_{gp}(A)$ is  finite.  See \cite{LM}, \cite{BB} and \cite{F}. \\ In the recent  paper \cite{BS} (to which we refer for an overview of the rich history of Question 1.1), Brion and Schr\"oer prove that any \textit{connected}   $G$ is isomorphic to $\AAut^0(X)$, for some projective, geometrically integral $F$-variety $X$.\\
This paper treats  the case of a \textit{linear} algebraic $F$-group $G$, possibly non-reduced.\\ The answer is then positive in full generality- see Theorem \ref{MainTh}. In the recent paper  \cite{Bg}, Bragg proves that every  finite étale $F$-group scheme is isomorphic to $\AAut(C)$, for $C$ a proper, smooth,  geometrically integral $F$-curve. In our work, it is unclear whether  assuming $G/F$ finite étale, could  lead to a simpler proof of Theorem \ref{MainTh}.\\
This paper is organised as follows. The main Theorem is stated in section \ref{SecThm}. Its proof occupies  Section \ref{SecProof}. 
Tools and intermediate results (most of which are unusual in positive characteristic) are developped in Sections 3 to 8.  Some of them are of independent interest. Here are two examples. \begin{enumerate}
    \item{Let  $X$ be the blow-up of a smooth $F$-variety $Y$, along a smooth closed subvariety $Z$.  Proposition \ref{PropAutBl} states that the infinitesimal automorphisms of $X$ are, as expected,  the infinitesimal automorphisms of $Y$ stabilising $Z$. } \item{Let $G$ be a linear algebraic group over  $F$. Proposition \ref{DividedTannaka} states that there exists a finite-dimensional $F$-vector space $W$, an integer $n \geq 1$, and a linear subspace $L \subset \Gamma^n_F(W)$, such that $G \simeq \Stab(L) \subset \PGL(W)$.  In \cite{M}, a related result is proved: $G$ is isomorphic to the stabiliser of a single   tensor of type $(2,1)$ (aka a non-associative finite-dimensional $F$-algebra). Can one use this,   instead of Proposition \ref{DividedTannaka}, to prove Theorem \ref{MainTh}? }\end{enumerate}
    
    The proof of Theorem \ref{MainTh} is considerably simpler when  $\carac(F)=0$, for two reasons.  First, divided powers may be  replaced by symmetric powers, and most algebraic results (e.g. Lemma \ref{LemFreeDiv}) become easy exercises. Second, whenever checking that a homomorphism $\phi: G \to H$ of linear algebraic $F$-groups is trivial (resp. injective, surjective), it suffices to prove that $\phi(\overline F): G(\overline F) \to H(\overline F)$ has the same property, as a morphism of abstract groups. Thus, differential calculus may be dismissed entirely. Sections  \ref{SecJet} and \ref{SectInfAut}  are not needed,   the length of the proof of  Proposition \ref{DividedTannaka} is halved,  and that of Lemma \ref{LemFreeDiv} is reduced tenfold.

\section{Statement of the Theorem.}\label{SecThm}
\begin{thm}\label{MainTh}
Let $G$ be a linear algebraic group over a field $F$.  \\ There exists a smooth projective $F$-variety $X$, such that $G$ is isomorphic to $\AAut(X)$, as a group scheme over $F$. More precisely, $X$ can be picked as the blow-up of a projective space, along a suitable smooth $F$-subvariety. 
\end{thm}

\section{Conventions, notation.}
Rings and algebras over them, are commutative with unit.\\ Denote by $F$ a field, with algebraic closure $\overline F$. Unless specified otherwise,  by ``$F$-vector space'' one means ``finite-dimensional $F$-vector space''.  Denote by $F[\epsilon]$, $\epsilon^2=0$, the $F$-algebra of dual numbers.  A variety over $F$ is a separated $F$-scheme of finite type. A linear algebraic group over $F$ is an affine $F$-variety, equipped with the structure of a group scheme over $F$. Equivalently, a linear algebraic group over $F$ is a closed $F$-sub-group scheme of $\GL_n$, for some $n \geq 1$.\\
Let $X$ be a variety over $F$. For an $F$-algebra $A$, denote by $X_A:=X \times_F A$ the $A$-scheme obtained from $X$ by extending scalars. Set $\overline X:=X\times_F \overline F$. \\The tangent sheaf $TX \to X$ is defined point-wise,  for every $F$-algebra $A$, by $$TX(A)=X(A[\epsilon]).$$ If $X$ is smooth over $F$, it is (the total space of) a vector bundle, dual to $\Omega^1(X/F)$.  \\A global section of the tangent sheaf is called a vector field on $X$. 
\subsection{Automorphism groups.}\label{AutGroups}For an  $F$-variety $X$, denote by $\AAut(X)$ the  automorphism group functor of $X$. For every $F$-algebra $A$, $\AAut(X)(A)$ is defined as the group of automorphisms of the $A$-scheme $X_A$. If $X/F$ is  projective, this functor is represented by a group scheme, locally of finite type over $F$. \\ For $X/F$ arbitrary, by  Lemma 3.1 of \cite{B}, there is a canonical isomorphism \[ H^0(X,TX) \stackrel \sim \to \Lie(\AAut(X)).\] 
Let $G/F$ be a group scheme, locally of finite type. If $G$ acts on the $F$-variety $X$, and for a closed subscheme $Z \subset X$, we use the notation  $\SStab_G(Z) \subset G$ for the  closed $F$-subgroup scheme defined by \[\SStab_G(Z) (A)=\{g \in G(A), g(Z_A)=Z_A\}, \]  for all  $F$-algebras $A$. That it is representable follows from \cite{DG}, II 1.3.6.
\subsection{Grassmannians.}

Let $V$ be an $F$-vector space. Pick an integer $d$, $0 \leq d \leq \dim(V)$. Denote by $\Gr(d,V)$ the Grassmannian of $d$-dimensional \textit{subspaces} of $V$. Set $\P(V)=\Gr(1,V)$. For $v \in V-\{0\}$, denote by $(v) \in \P(V)(F)$ (or abusively $v$ if no confusion arises) the line directed by $v$. Recall that, for $E \in \Gr(d,V)(F)$, the tangent space $T_E(\Gr(d,V))$ is naturally isomorphic to $\Hom_F(E, V/E)$.
\subsection{Weil restriction.}\label{WeilRes}
Recall the following important tool.
\begin{defi}
   Let $A$ be a finite $F$-algebra. Let $Y$ be a quasi-projective scheme over $A$. Denote by $R_{A/F}(Y)$ the Weil restriction of $Y$. It is a  quasi-projective variety over $F$, characterised by the formula, for every $F$-algebra $B$:  \[ R_{A/F}(Y)(B)=Y(B \otimes_F A).\] 
\end{defi}

\subsection{Symmetric and divided powers.} 
Let $V$ be a vector space over $F$. Define \[ V^\vee :=\Hom_F(V,F).\] For each $n \geq 1$, define \[ \Sym^n(V):=H_0(S_n,V^{\otimes n})\] and \[ \Gamma^n(V):=H^0(S_n,V^{\otimes n}),\] where coinvariants and invariants are taken w.r.t. the natural permutation action of the symmetric group $S_n$ on  $V^{\otimes n}$. These are, respectively, the $n$-th symmetric power and the $n$-th divided power of $V$. For $v \in V$, set \[ [v]_n:=v \otimes v \otimes \ldots \otimes v \in \Gamma^n(V).\] These are called pure symbols. If $\vert F\vert \geq n$, they span  $\Gamma^n(V)$.\\  There is a canonical non-degenerate pairing of $F$-vector spaces $$ \Gamma^n(V^\vee) \times \Sym^n(V) \to F,$$ $$ ([\phi]_n, x_1 x_2 \ldots x_n) \mapsto \phi(x_1) \ldots \phi(x_n). $$ It is perfect- whence an isomorphism \[\Gamma^n(V^\vee) \stackrel \sim \to  \Sym^n(V)^\vee.\] 

Denote by $\Sym(V)=\bigoplus_{n \in \N} \Sym^n(V)$ and $\Gamma(V)=\bigoplus_{n \in \N} \Gamma^n(V)$ the symmetric and divided powers algebras of $V$, respectively. Relations in $\Gamma(V)$, arise from the motto  ``$[v]_n= \frac {v^n} {n!}$''. These are:
\begin{enumerate}
    \item{$[v]_0=1$,} \item{$ [v+v']_n= \sum_0^n [v]_i [v']_{n-i}$,} \item{$ [\lambda v]_n=\lambda^n[v]_n$,} \item{ $[v]_n[v]_m= {n+m \choose n}[v]_{n+m}$.}
\end{enumerate}For details,  see \cite{Ro}.\\Let $(e_1, \ldots, e_d)$ be an $F$-basis of $V$. Then,   $\Gamma^n(V)$ inherits a (canonical) basis, consisting of symbols $[e_1]_{a_1} \ldots [e_d]_{a_d}$, where $a_i \geq 0$ and $a_1+ \ldots + a_d=n$. Dually,  $\Sym^n(V)$ inherits its usual monomial basis, consisting of tensors $e_1^{a_1} \ldots e_d^{a_d}$. \\There are two natural arrows \[\Sym^n(V) \to  \Gamma^n(V),\] \[ v_1 v_2 \ldots v_n \mapsto [v_1]_1 [v_2]_1 \ldots [v_n]_1,\] and \[\Gamma^n(V) \to  \Sym^n(V),\] \[[v]_n \mapsto v^n.\] Their composites equal $n! \Id$. Hence, if $\carac (F)=0$ or $p>n$, they are isomorphisms.

\section{Jet spaces via infinitesimal Weil restriction.}\label{SecJet}

 ``One-dimensional'' jet spaces (i.e. with values in $F[X]/X^n$, for some $n \geq 1$) are a famous tool in many branches of geometry. However, they would not suffice here. In this section,  we offer  a self-contained exposition of what we shall actually need.  
 \begin{defi}
     Denote by $A$   a finite local $F$-algebra with residue field $F$, by  $\mathcal M \subset A$ its maximal ideal,  and  by $\rho: A \to A/\mathcal M=F$ its residue homomorphism, which is a retraction of the inclusion $F \hookrightarrow A$.  
 \end{defi}

\begin{defi}(Jet spaces.)\label{DefiJet} \\
let $q: A \to A'$ be a homomorphism of finite local $F$-algebras with residue field $F$. For any $F$-variety $X$, $q$ induces a morphism of $F$-varieties $$q_*:  R_{A/F}(X) \to R_{A'/F}(X).$$  Formula: for an  $F$-algebra $B$ and  for $$ x \in X(B \otimes_F A),$$ $q_*(x)$ is defined as $$(q \circ x) \in X(B \otimes_F A').$$\\ 
Let $G$ be a contravariant group functor, on affine $F$-varieties.\\  Define a group functor   $\mathbf{J}(G,\rho)$ by $$\mathbf{J}(G,\rho)(B):=\Ker(G(A \otimes_F B) \stackrel {\rho_*} \to G(B)),$$ for every $F$-algebra $B$. For an $F$-variety $X$, set $$\mathbf{J}(X,\rho):=\mathbf{J}(\AAut(X),\rho).$$ 
\end{defi}

\begin{ex}
    If $A=F[\epsilon]$, then  $\rho_*$  is  the tangent sheaf $TX \to X$, and $\mathbf{J}(X,\rho)=H^0(X,TX)=\Lie(\AAut(X)).$
\end{ex}
\begin{lem}\label{LemWeilRes}
 Keep notation and assumptions above.    Assume moreover, that $X$ is a smooth $F$-variety. 
 Consider a diagram  \[ \xymatrix {A_1 \ar[r]^q \ar[d]^{\rho_1} & A_2 \ar[d]^{\rho_2} \\ F \ar@{=}[r] & F,}\] of epimorphisms of  finite local $F$-algebras with residue field $F$. \\ Denote by $\mathcal M_i$ the maximal ideal of $A_i$, and set $\mathcal I:=\Ker(q)$.\\ Assume that $\mathcal I \mathcal M_1=0$. Then, the morphism  of $F$-varieties $$q_*:   R_{A_1/F}(X) \to  R_{A_2/F}(X)$$  is   a torsor under the (pull-back via $\rho_2^*$ of the) vector bundle  $TX \otimes_F \mathcal I$.\\
 Thus, the morphism $\rho_*$ is a composite of torsors under the vector bundle $TX$.\\As such, it is affine and smooth.

\end{lem}

\begin{dem}

 The assertion is local on the smooth $F$-variety $X$, so that one may assume $X=\Spec(R)$ affine. Then  $R_{A_i/F}(X)$ ($i=1,2$) is affine as well. Let $B$ be an $F$-algebra. By the infinitesimal lifting criterion for smooth morphisms (\cite{SP},  tag 37.11.7), one sees that the map $$ q_*(B):  X(A_1 \otimes_F B) \to   X(A_2 \otimes_F B) $$ is onto. Let $x_1, y_1 \in  X(A_1 \otimes_F B)$ be such that $q_*(B)(x_1)=q_*(B)(x_2)$.  Set $$x_0:= (\rho_1)_*(B)(x_1)=(\rho_2)_*(B)(x_2)\in X(B).$$Consider $x_1,y_1$ (resp. $x_0$) as homomorphisms of $F$-algebras $R \to A_1 \otimes_F B$ (resp. $R \to B$), and form the difference $$\delta:= (y_1-x_1): R \to  A_1 \otimes_F B.$$  This is a priori just an $F$-linear map. Since $q_*(B)(x_1)=q_*(B)(x_2)$, it takes values in $ \mathcal I \otimes_F B \subset  A_1 \otimes_F B$.  Consider $\mathcal I  \otimes_F B$ as an $R$-module via $x_0$, treating $\mathcal I$ just as an $F$-vector space. One then checks that $$\delta:  R \to   \mathcal I  \otimes_F B$$ is an $F$-derivation. Conversely,  assume given  a homomorphism of $F$-algebras $$x_1: R \to A_1 \otimes_F B.$$ Denote its reduction mod $\mathcal M$ by $x_0: R \to  B$. Pick   a derivation $$\delta:  R \to   \mathcal I  \otimes_F B,$$ where the target is considered as an $R$-module as above. Then $$y_1:=(x_1+\delta):  R \to A_1 \otimes_F B$$  is a homomorphism of $F$-algebras, such that  $q_*(B)(x_1)=q_*(B)(x_2)$. \\This completes the description of the torsor structure.\\  For the last assertion, one may assume $\mathcal M \neq 0$. Then, the minimal non-zero ideals of $A$ are one-dimensional, and killed by $\mathcal M$. Thus,  $\rho$ can be written as a composite  of epimorphisms of finite local $F$-algebras, $$A =A_n \to A_{n-1} \to \ldots \to A_2 \to A_1=F ,$$ where  $\dim_F(A_i)=i$. Induction on $\dim_F(A)$ then applies.
\end{dem}

\begin{rem}
  Let $V$ be an $F$-vector space.  Consider its affine space $$X=\A_F(V):=\Spec(\Sym^*(V^\vee)).$$ Treating $\rho$ as a linear form on the $F$-vector space $A$, $\rho_*$ is simply $$\A_F(V \otimes_F A) \to \A_F(V),$$ $$w \mapsto (\Id_V \otimes \rho)(w).$$  Thus, it is a trivial fibration in affine spaces.  
\end{rem}
\begin{rem}

If $X$ is affine, using  the preceding Lemma, and vanishing of coherent cohomology over an affine base, one sees that $\rho_*$ is a trivial $\A^N$-fibration, as well.
\end{rem}

\begin{prop}\label{PropJet}
    Let $X$ be an $F$-variety. 
    There is a functorial isomorphism between $\mathbf{J}(X,\rho)$, and the functor of sections of the morphism of $F$-schemes $R_{A/F}(X_A) \stackrel {\rho_*} \to X. $ 
\end{prop}

\begin{dem}
Let us describe, for every $F$-algebra $B$, a functorial bijection$$\mathbf{J}(X,\rho)(B) \stackrel \sim \to  \{s: X_B \to R_{A\otimes_F B /B}(X_{A\otimes_F B}), \; \rho_* \circ s =\Id\}.$$
For simplicity,   we assume $B=F$; the construction actually works in general.\\ Giving a section $s: X \to R_{A/F}(X_A)$ amounts to giving a morphism of $A$-schemes $f: X \times_F A \to X \times_F A$.  Assuming that $\rho_*(f): X \to X$  is the identity, one then just needs to show that $f$ is  an iso. Since $\mathcal M$ is nilpotent, one sees that, as a homeomorphism of the topological space $X \times_F A$, $f$ is the identity. Let $(U_i)$ be a covering of $X$ by open affines. From what was just said, $f$ restricts to morphisms of $A$-schemes $f_i: U_i \times_F A \to U_i\times_F A$. By a straightforward glueing argument, one thus reduces to the case $X$  affine. One may then use Lemma \ref{LemNakArt} below, applied to the homomorphism of $A$-algebras $$\Phi: \mathcal O_X(X) \otimes_F A  \to \mathcal O_X(X) \otimes_F A,$$ which is such that $f=\Spec(\Phi)$. Note that $\Phi$ is regarded here as a morphism between free $A$-modules. Since $\phi=\Id$ is an iso, one concludes that $\Phi$ is an iso. Hence  $f$ is an iso of $A$-schemes, as desired. 
\end{dem}

\begin{coro}\label{CoroJAlg}
    Assume that $X$ is a smooth projective $F$-variety. Denote by $\AAut(X)^0 \subset \AAut(X)$ the connected component of the identity. It  is  a group scheme,  of finite type over $F$. Then $\mathbf{J}(X,\rho)= \mathbf{J}(\AAut(X)^0,\rho)$ is a smooth, connected  and unipotent linear algebraic $F$-group. Moreover, it  is  $F$-split. In other words, it has a composition series with quotients $\G_a$.
\end{coro}

\begin{dem}
    This follows from   Lemma \ref{LemWeilRes} and Proposition \ref{PropJet}.
\end{dem}

\begin{lem}\label{LemExtJ}
    Let $\iota: Z \hookrightarrow Y$ be a closed immersion of smooth affine $F$-varieties. \\
  Denote by  $\mathbf J(\iota,\rho)\subset \mathbf J(Y,\rho)$ (resp. $\mathbf J_0(\iota,\rho)\subset \mathbf J(Y,\rho)$) the sub-group  functor formed by infinitesimal automorphisms $f$, such that $f_{\vert Z}$ factors through $\iota$ (resp. $f_{\vert Z}=\iota$). There is an exact sequence (of group functors) $$ 1 \to  \mathbf J_0(\iota,\rho)\to  \mathbf J(\iota,\rho)\xrightarrow{f \mapsto f_{\vert Z}}  \mathbf J(Z,\rho) \to 1.$$

\end{lem}

 \begin{dem}
    Let us show  that  $$ 1 \to  \mathbf J_0(\iota,\rho)(F)\to  \mathbf J(\iota,\rho)(F)\xrightarrow{f \mapsto f_{\vert Z}}  \mathbf J(Z,\rho)(F) \to 1$$ is exact, as a sequence of abstract groups. The  same proof  works to show exactness for points in an arbitrary $F$-algebra $R$.  The only non-trivial part is surjectivity, which we check by  induction on $\dim_F(A),$ using  Lemma \ref{LemWeilRes}, of which we adopt notation, and the description of $\mathbf J(X,\rho)$ provided in Proposition \ref{PropJet}. Let $f_1: Z \to R_{A_1 /F}(Z)$. By induction, $f_2:=q_* \circ f_1$ extends, to $\tilde f_2: Y \to R_{A_2 /F}(Y)$. Since  $Y$ is affine, $H^1(Y, TY \otimes_F \mathcal I)=0,$ so that $\tilde f_2$ lifts via $q_*$, to $g_1: Y \to R_{A_1 /F}(Y)$. Consider $(g_1)_{\vert Z}: Z \to R_{A_1 /F}(Y)$. Via $q_*$, it is sent to $f_2$. Thus, there exists a unique $\epsilon \in H^0(Z,TY \otimes_F \mathcal I),$ such that $(g_1)_{\vert Z}+\epsilon=f_1.$ Again, since $Y$ is affine, $\epsilon$ extends, to  $\tilde \epsilon \in H^0(Y,TY \otimes_F \mathcal I).$ Then, $\tilde f_1:=g_1+ \tilde \epsilon$ is the sought-for extension of $f_1$.\\
    
   \end{dem}

The following result is standard. Lacking a reference, a proof is included.
\begin{lem}(Improved Nakayama's, for Artinian rings.)\label{LemNakArt}\\
 Let $M,N$ be $A$-modules, and let $\Phi:M \to N$ be an $A$-linear map. Denote by  $\phi:M\otimes_A F \to N\otimes_A F$ the induced $F$-linear map. The following holds. 
    
    \begin{enumerate}

        \item{If $\phi$ is surjective, so is $\Phi$. }

        \item{Assume that $N$ is  a  free $A$-module,  possibly of infinite rank. \\Then if $\phi$ is injective, so is $\Phi$.} 
    \end{enumerate}
\end{lem}

\begin{dem}
  To prove (1), proceed  by induction  on the smallest $k \geq 1$, such that $\mathcal M^k M=\mathcal M^k N=0$. Case $k=1$ is clear. Assume that $\mathcal M^{k+1}M=\mathcal M^{k+1}N=0$, and that $\phi$ is onto. By induction, $(\phi/\mathcal M^k): M/ \mathcal M^k M \to N /  \mathcal M^k N$ and $\phi_{\vert \mathcal M^k M}:  \mathcal M^k M \to \mathcal M^k N$ are onto. By  dévissage, one then sees that $\Phi$ is onto.\\Let's prove  (2) by induction on the length of $A$. It suffices to show the following. \\ Let $ a \neq 0
\in \mathcal M$ be such that $a \mathcal M=0$. Assume that  the $(A/aA)$-linear map $$\Phi_a:=(\Phi/a): M/aM \to N/aN$$ is injective (in addition to injectivity of $\phi$). Then $\Phi$ is injective. \\To prove this assertion, pick $x \in \Ker(\Phi)$. Since $\Phi_a$ is injective, one gets $x \in aM$. Write $x=am$. Since $a \Phi(m)=\Phi(x)=0 \in N$,  since $N$ is a free $A$-module, and since $A-\mathcal M=A^\times$, it must be the case that $\Phi(m) \in \mathcal M N$. Hence $\phi(m)=0$. Injectivity of $\phi$ then implies $m \in \mathcal M M$. Thus, $x \in a \mathcal M M= 0$.
\end{dem}

\section{Infinitesimal automorphisms of blow-ups.}\label{SectInfAut}

The following  improves on \cite{F}, Lemma 4.2. 
\begin{prop}\label{PropAutBl}
    Let $\iota: Z \hookrightarrow Y$ be a closed immersion of smooth $F$-varieties. Denote by $$\beta: X:=\Bl_Z(Y) \to Y$$ the blow-up of $Y$ along $Z$, and by $$i: E \hookrightarrow X$$ the exceptional divisor.  Let $A$ be an $F$-algebra. There is a natural monomorphism  $$\phi: \Stab_{\AAut(Y)}(Z)(A) \to \AAut(X)(A) .$$ Assume  that each irreducible component of $Z$ has codimension $\geq 2$ in $Y$. \\
 Let $A$ be a finite local $F$-algebra, with residue homomorphism $\rho: A \to F$.\\ Then $\phi$ induces an iso $$\Phi: \mathbf{J}(\Stab_{\AAut(Y)}(Z),\rho) \stackrel \sim \to  \mathbf{J}(X,\rho).$$
   
\end{prop}

\begin{dem}
 Recall that formation of blow-ups is functorial and commutes to base-change. Precisely,  let $f: Y\times_F A \to  Y\times_F A $ be an automorphism of $A$-scheme, preserving the subscheme $Z \times_F A$. By functoriality of the blow-up,  $f$ induces an $A$-automorphism of $X \times_F A$. This  provides the definition of $\phi$.  Assume that $f \in \Stab_{\AAut(Y)}(Z)(A)$ is such that $\phi(f)=\Id$. Consider the commutative diagram $$\xymatrix{X_A \ar[r]^{\beta_A} \ar[d]^{\Id} & Y_A \ar[d]^f\\ X_A \ar[r]^{\beta_A}  & Y_A. }$$  Since  $\beta_A$ is surjective,  one sees that $f$, as  a continuous map, is the identity.  Checking that $f=\Id$ becomes local  on $Y$, so that one may assume  $Y=\Spec(B)$ is affine and connected. Then $B$ is integral, because $Y$ is smooth over $F$. Denote by $I\subset B$ the ideal defining $Z$. Then $X:=\Proj (R)$, where $R:=\bigoplus_{n=0}^\infty I^n$. It is covered by the open affines $\Spec(R[\frac 1 f]_0)$, for $0 \neq f \in I$. Since $B$ is integral, the natural arrow $B \otimes_F A \to  R[\frac 1 f]_0 \otimes_F A$  is injective. The claim follows.  \\For the second assertion,  note that elements of the source and target of $\Phi$ are topologically the identity. Thus, the question is local on $Y$, so that one may assume $Y$ (and hence $Z$) affine. We use (and adopt notation of) Proposition \ref{PropJet}. Note  that  $\mathbf{J}(\Stab_{\AAut(Y)}(Z),\rho) \subset \mathbf J (Y,\rho) $ is the sub-functor  $\mathbf J (\iota,\rho)$ of Lemma \ref{LemExtJ}.  By \cite{F}, Lemma 4.2 (or more accurately, its proof), one knows that
 $$ H^0(X,TX) =\Ker(H^0(Y,TY) \to H^0(Z, N_{Z/Y})).$$ Equivalently, a vector field on $X$ is the same thing as a vector field on $Y$, restricting to a vector field on $Z$.  Via this identification, the torsor structures on both sides are  automatically compatible with $\Phi$. One can then proceed by induction on $\dim_F(A)$ again, using  Lemma \ref{LemWeilRes}, which provides an exact sequence $$ 0 \to  H^0(Y, TY \otimes_F \mathcal I) \otimes_F R \to \mathbf J(Y, \rho_1)(R) \to \mathbf J(Y, \rho_2)(R) \to 1, $$  functorial in the $F$-algebra $R$.  For simplicity, let us work with $F$-points- the case of $R$-points being the same. Consider the natural diagram (of abstract groups) $$\xymatrix{0 \ar[r] & \Ker(H^0(Y,TY) \to H^0(Z, N_{Z/Y}))\otimes_F \mathcal I \ar[r]  \ar[d]^{\simeq} & \mathbf J(\iota,\rho_1)(F)\ar[r]  \ar[d] & \mathbf J(\iota,\rho_2)(F)  \ar[r] \ar[d]^{\simeq} &  1\\0 \ar[r] & H^0(X,TX )\otimes_F \mathcal I \ar[r] & \mathbf J(X,\rho_1)(F)\ar[r]  & \mathbf J(X,\rho_2) (F) \ar[r] &  1 .  } $$ [Note that the vertical arrows are given by functoriality of the blow-up. In the bottom line, we used $H^0(X,TX \otimes_F \mathcal I )=H^0(X,TX )\otimes_F \mathcal I$, and a similar fact in the upper line, which hold   because $\dim_F(\mathcal I) < \infty$.] In this diagram, lines are clearly exact, except possibly on their right. To conclude, it remains to prove surjectivity of $\mathbf J(\iota,\rho_1)(F) \to \mathbf J(\iota,\rho_2)(F)$.  Using the exact sequence of Lemma \ref{LemExtJ},  a diagram chase reduces this to checking  surjectivity of $\mathbf J_0(\iota,\rho_1) (F)\to \mathbf J_0(\iota,\rho_2)(F) $. Pick $f_2 \in  \mathbf J_0(\iota,\rho_2)(F)  \subset \mathbf J(Y,\rho_2)(F)$. Extend it (via $q_*$) to  $g_1 \in  \mathbf J(Y,\rho_1)(F).$  Then $q_*((g_1)_{\vert Z})=\iota$, so that $(g_1)_{\vert Z}=\iota+\epsilon$, for $\epsilon \in H^0(Z, TY \otimes_F \mathcal I)$. Since $Y$ is affine, $\epsilon$ extends, to $\tilde \epsilon \in H^0(Y, TY \otimes_F \mathcal I)$. Then, $f_1:=g_1- \tilde \epsilon$ is the sought-for lift of $f_2$.

\end{dem}

\section{Automorphisms of blow-ups of projective space, via Chow rings.}
We begin with gathering, from \cite{Fu},  material on blow-ups and their Chow rings.
\begin{prop}\label{ChowBlow}
    Let $\iota: Y \hookrightarrow Z$ be a closed immersion between smooth geometrically integral $F$-varieties,  of codimension $c \geq 2$. Denote by $$\beta: X:=\Bl_Y(Z) \to Z$$ the blow-up of $Z$ along $Y$, and by $e: E \hookrightarrow X$ the exceptional divisor. 
    
    \begin{enumerate}

    \item{The restriction $\beta_{\vert X-E}: (X-E) \to (Z-Y)$ is an isomorphism, providing a natural arrow $$\phi: \Pic(X) \to  \Pic(Z-Y)=\Pic(Z).$$ $$L \mapsto L_{\vert X-E}. $$}
        \item{The morphism $$\pi:=\beta_{\vert E}: E \to Y$$ is the projective bundle of the  normal bundle $N_{Y/Z}$.\\ Denote by $\mathcal O_E(1)$ its twisting sheaf, and set  $ \zeta:=c_1(\mathcal O_E(1)) \in \CH^1(E).$\\ The normal bundle $N_{E/X}$ is canonically isomorphic to $\mathcal O_E(-1)$.\\ For all $i \geq 1$,$$ [E]^i=(-1)^{i-1} e_*( \zeta^{i-1}) \in \CH^i(X).$$ }

        \item{(Projective bundle formula for $\pi$).  The arrow $$ \bigoplus^{i}_{j=i+1-c}\CH^j(Y)  \to \CH^i(E),$$  $$(x_j) \mapsto \sum \pi^*(x_j).\zeta ^{i-j} $$ is an iso, for every $i \geq 1$. \\ In particular, for $i=1$, the natural arrow $$ \Pic(Y) \bigoplus \Z \to \Pic(E),$$  $$ (L,a) \mapsto \pi^*(L)+ \mathcal O(a) $$ is an iso. Denote the  projection on the second factor by $$\Pic(E) \to \Z, $$ $$ L \mapsto d(L).$$}
    \item{ The natural arrow $$ \Pic(Z) \bigoplus \Z \to \Pic(X),$$ $$ (L, a) \mapsto \beta^*(L)+\mathcal O_X(aE)$$ is an iso. Its inverse is given by $$\Pic(X) \to \Pic(Z) \bigoplus \Z, $$ $$ L \mapsto (\phi(L) , -d(e^*(L)).$$}
\item{More generally, for $i \geq 1$ there is a natural exact sequence   $$0 \to \CH^{i-c}(Y) \to \CH^i(Z) \bigoplus \CH^{i-1}(E) \xrightarrow{\sigma}\CH^i(X) \to 0,$$ with injection  given by $$ w \mapsto (-\iota_*(w),\pi^*(w).\zeta^{c-1}), $$ and surjection given by $$(u,v) \mapsto \beta^*(u)+ e_*(v). $$If $i < c$, this boils down to an isomorphism   $$ \CH^i(Z) \bigoplus \CH^{i-1}(E) \stackrel \sim \to \CH^i(X).$$} \item{Let $ i \geq 1$, and   $u \in \CH^i(Z)$. Via item (5), the product $$\beta^*(u).[E]  \in \CH^{i+1}(X) $$  equals $\sigma(0,\pi^*(\iota^*(u))).$}
    \end{enumerate}

\end{prop}
\begin{dem}

Items (1), and  the first two assertions of (2),  are  standard features of blow-ups. The self-intersection formula  for $[E]$ (\cite{Fu}, Section 6.3), and the projection formula (\cite{Fu}, Example 8.1.1), then prove the last formula of (2) by induction on $i$: $$ [E]^{i+1}= [E].[E]^i=(-1)^{i-1}[E].e_*(\zeta^i)= (-1)^{i-1} e_*(e^*([E]).\zeta^i)=(-1)^i e_*(\zeta^{i+1}).$$    [The starting case $i=1$ holds by definition.]\\ Item (3) is \cite{Fu}, Theorem 3.3.b. Item (5) is  Proposition 6.7.e of \cite{Fu} (note the explicit formulas in its proof). Item (4) is a particular case of (5), for $i=1$.\\
   Observe that $\beta \circ e=\iota \circ \pi. $ Item (6) follows, using   the projection formula:$$ \beta^*(u).[E] =\beta^*(u). e_*(1_E)=e_* (e^*(\beta^*(u)))=e_*(\pi^*(\iota^*(u))). $$
\end{dem}

The content of the following two  Propositions is that, under suitable  assumptions, the automorphism group of a blow-up in projective space, is ``as  naively expected.''

\begin{prop}\label{PhiIso}
     For $N \geq 6$, let  $Y_1,Y_2 \subset \P^N$ be  disjoint  smooth  closed $F$-subvarieties, geometrically integral and of  dimensions in $ [1,N-3]$. \\Denote by $$\beta: X:=\Bl_Y(\P^N) \to \P^N$$ the blow-up of $\P^N$ along $Y:=Y_1 \coprod Y_2$, and by $e_i: E_i \hookrightarrow X$ the exceptional divisor lying above $Y_i$, $i=1,2$.  By functoriality of the blow-up, there is a natural  homomorphism of abstract groups $$\Phi: \Stab_{\Aut(\P^N)}(Y_1) (F)\cap   \Stab_{\Aut(\P^N)}(Y_2)(F)  \to \Aut(X)(F). $$ If the $F$-varieties $E_1$ and $E_2$ are not isomorphic, then $\Phi$  is an isomorphism.
\end{prop}
\begin{dem}
That $\Phi$ is injective is straightforward. Let us check surjectivity.
  Let $f: X \to X$ be an $F$-automorphism. Observe that $\beta$ is the composite $$ X \xrightarrow{\beta_2} X_1 \xrightarrow{\beta_1} \P^N,  $$ where $\beta_1$ (resp. $\beta_2$)  is the blow-up of $\P^N$ along $Y_1$ (resp. of $X_1$ along $\beta_1^{-1}(Y_2)$). By item (4) of Proposition \ref{ChowBlow}, applied two times, one gets that $\CH^1(X)$ is a free $\Z$-module of rank $3$, with basis ($\beta^*([H]), [E_1]$,$[E_2]$), where $H \subset \P^N$ is a hyperplane. Since $c \geq 3$, item  (5) (applied two times, to $\beta_1$ and $\beta_2$) provides a natural iso  $$\CH^2(X) \simeq \Z.[\beta^*(H)]^2  \bigoplus \CH^1(E_1) \bigoplus \CH^1(E_2).$$   Assume that $f(E_i) \neq E_j$, for all $\{i,j \} \subset \{1,2\}$. Since $f(E_i) \subset X$ is an effective divisor not contained in $E_1 \coprod E_2$,  the last formula of item (4) then yields a decomposition, for $i=1,2$, $$[f(E_i)] =a_i[\beta^*(H)]-b_i[E_1]-c_i[E_2] \in \CH^1(X),$$  with $a_i \geq 1$,  and $b_i, c_i \geq 0$. Using  item (6) two times  (exchanging the roles of $E_1$ and $E_2$), one gets, for $i=1,2$,  $$ [E_1]. [E_2] \in \CH^1(E_i) \subset \CH^2(X),$$ w.r.t. the direct sum decomposition above. Thus   $ [E_1]. [E_2]=0$. One also computes  $$[f(E_1)]. [f(E_2)] = (a_1[\beta^*(H)]-b_1[E_1]-c_1[E_2]). (a_2[\beta^*(H)]-b_2[E_1]-c_2[E_2])$$ $$=(a_1a_2[\beta^*(H)]^2,\ast,\ast),$$ where it is needless to know the expression of the second and third components. It suffices to observe that  $[f(E_1)]. [f(E_2)] \neq 0$, whereas $[E_1]. [E_2]=0$. This is impossible, since $f$ induces a ring automorphism of  $\CH^*(X)$. Consequently, it must be the case that $f(E_i)=E_j$ for some $\{i,j \} \subset \{1,2\}$. Then $i=j$, by the assumption made on $E_1$ and $E_2$. Say $i=j=2$, so that  $f(E_2) =E_2$.  Assume that $f(E_1) \neq E_1$. Then, as above,  one may write  $$[f(E_1)] =a_1[\beta^*(H)]-b_1[E_1]-c_1[E_2] \in \CH^1(X),$$  with $a_1 \geq 1$,  and $b_1, c_1 \geq 0$. Compute: $$[f(E_1)].[f(E_2)]=(a_1[\beta^*(H)]-b_1[E_1]-c_1[E_2]).[E_2]$$ $$ =a_1[\beta^*(H)].[E_2]-c_1[E_2].[E_2] \in \Pic(E_2) \subset \CH^2(X), $$ w.r.t. the direct sum decomposition above. Via the projection formula for the projective bundle $\pi_2: E_2  \to Y_2$ (items (2) and (3) of Proposition \ref{ChowBlow}),  one gets $$a_1[\beta^*(H)].[E_2]-c_1[E_2].[E_2]=(\iota_2^*(a_1[H]),c_1)\in (\Pic(Y_2) \bigoplus \Z)\simeq \Pic(E_2),$$ where $\iota_2 : Y_2 \hookrightarrow \P^N$ is the closed immersion, and using
  item (6) with $u:=[H]$. Since $a_1 \geq 1$, the divisor class  $\iota_2^*(a_1 [H]) \in \Pic(Y_2)$ is ample, on the positive-dimensional projective variety $Y_2$, hence is non-zero. It follows that  $[f(E_1)].[f(E_2)] \neq 0$, contradiction. Thus, $f(E_1) = E_1$ and $f(E_2) = E_2$. Then, $f$ restricts to an automorphism $g$ of $X-E_1-E_2$, which by item (1) is an open subvariety of $\P^N$, with complement $Y$ of codimension $\geq 2$.  By Lemma \ref{LemAutPN}, $g$ indeed extends to an automorphism of $\P^N$, which necessarily fixes $Y_1$ and $Y_2$ separately.

\end{dem}
\begin{rem}
    Under the same assumptions, Proposition \ref{PhiIso} can be generalised to a blow-up of any number of disjoint smooth subvarieties.
\end{rem}

  \begin{lem}\label{LemAutPN}
     Let $Y \subset \P^N$ be a closed $F$-subvariety, of codimension $\geq 2$. Set $U:=\P^N-Y$.  Every $F$-automorphism of $U$ extends to an automorphism of $\P^N$.
  \end{lem}

  \begin{dem}
   Recall that, on a normal $F$-variety, regular functions are invariant upon removing a closed subvariety of codimension $\geq 2$. The same property then holds for global sections of line bundles, and one  also infers that $\Pic(U)=\Pic(\P^N)=\Z.[\mathcal O(1)]$.  One can then reproduce the classical proof that $\Aut(\P^N)(F)=\PGL_{N+1}(F)$, with  $U$ in place of $\P^N$. Here are details. Let $g$ be an $F$-automorphism of $U$. Then $g^*([\mathcal O(1)])$ is ample and generates  $\Pic(U)$;  hence   $g^*(\mathcal O(1))\simeq \mathcal O(1)$. Fix such an iso of line bundles, and consider the effect of $g^*$ on $$H^0(U,\mathcal O(1))=<X_0, \ldots, X_N>=F^{N+1}.$$ This gives a well-defined $\tilde g \in \PGL_{N+1}(F)$- the sought-for extension of $g$.
  \end{dem}
The following instructive exercise concludes this section.  The proof given is by counting points over finite fields, which is more elementary than by Chow groups.
\begin{lem}\label{LemEnd}
 For $i=1,2$, let $a_i,m_i \geq 2$ be integers, and let $V_i$  be a vector bundle of rank $m_i$  over $\P^{a_i-1}_F$. Denote by $\P(V_i) \to \P^{a_i-1}_F$ the corresponding projective bundles.  Assume that $a_1 \neq a_2$, and that $\P(V_1)$ and $\P(V_2)$  are isomorphic as  $F$-varieties. Then $m_1=a_2$ and $m_2=a_1$.
 \end{lem}

 \begin{dem}
     By a classical ``spreading out'' argument, one may assume that $F=\F_q$ is a finite field. Indeed, there exists a sub-ring $R \subset F$, which is a $\Z$-algebra of finite type, such that all data in the Lemma are defined over $R$. More precisely, for $i=1,2$ there is  a vector bundle $\mathcal V_i$ of rank $m_i$  over $\P^{a_i-1}_R$, such that the projective bundles $\P(\mathcal V_1)$ and $\P(\mathcal V_2)$  are isomorphic as $R$-schemes. Specialising at a closed point of $\Spec(R)$, one gets a similar data over a finite field, as claimed. Consider the morphism of $\F_q$-varieties $\P(V_i) \to \P^{a_i-1}_{\F_q}$. It induces a surjection of finite sets $$\P(V_i)(\F_q) \to \P^{a_i-1}_{\F_q}(\F_q),$$ with fibers (non-canonically isomorphic to)  $\P_{\F_q}^{m_i-1}$.  Counting points, one gets $$\P(V_i)(\F_q)=\frac {(q^{a_i}-1)(q^{m_i}-1)}  {(q-1)^2}.$$ Since the $\F_q$-varieties $\P(V_1)$ and $\P(V_2)$  are isomorphic, one has  $$\frac {(q^{a_1}-1)(q^{m_1}-1)}  {(q-1)^2}=\frac {(q^{a_2}-1)(q^{m_2}-1)}  {(q-1)^2}.$$  For  $n \geq 1$, extend scalars to $\F_{q^n}$    to get the same formula, with $q^n$ in place of $q$. \\ Thus,  $$\frac {(X^{a_1}-1)(X^{m_1}-1)}  {(X-1)^2}=\frac {(X^{a_2}-1)(X^{m_2}-1)}  {(X-1)^2} \in \Q(X),$$ and the conclusion follows.
 \end{dem}

\section{Divided powers to the rescue of projective geometry.}

One can think of the results this section,  as a characteristic-free version of polarity. If $\carac(F)=0$, many of these boil down to  facts found in \cite{D}, chapter 1.

\subsection{Veronese embedding.}

Here is a convenient coordinate-free definition of the Veronese embedding. Up to the choice of a basis, it agrees with the usual one.
\begin{defi}
Let $V$ be an $F$-vector space. Let $n \geq 1$ be an integer.
The arrow of $F$-varieties \[\Ver_n: \P(V) \to \P(\Gamma^n(V)),\] \[ v  \mapsto [v]_n.\] is a closed embedding, called  the $n$-th Veronese embedding.
\end{defi}

\begin{prop}\label{AutVeronese}
Let $V$ be a finite dimensional $F$-vector space. Let $n \geq 1$ be an integer.
Consider the $n$-th Veronese embedding \[ \P(V) \xrightarrow{\Ver_n} \P(\Gamma^n(V)).\]
The natural arrow  \[ \AAut(\P(V))=\PGL(V) \to  \PGL(\Gamma^n(V))=\AAut(\P(\Gamma^n(V)))\] induces an isomorphism of linear algebraic $F$-groups \[\phi: \PGL(V) \to  \Stab_{\PGL(\Gamma^n(V))}(\P(V) \subset \P(\Gamma^n(V))).\]

\end{prop}
\begin{dem}
Can assume  $\dim(V) \geq 2$, and $F$ infinite.
Let $A$ be an $F$-algebra. Pick \[ f \in \Stab_{\AAut(\P(\Gamma^n(V)))}(\P(V))(A).\] Then $f$ restricts to an automorphism of the $A$-scheme $\P(V) \times_F A$; that is, to an element $\psi(f) \in \PGL(V)(A)$. By Yoneda's Lemma, this defines an $F$-morphism $$\psi: \Stab_{\AAut(\P(\Gamma^n(V)))}(\P(V)) \to \PGL(V),$$ which is a retraction of $\phi$. Hence, $\phi$ is an embedding. Since its source is smooth, it  suffices to show that every element  \[ f \in \Ker(\psi)(F[\epsilon])\]  lies in the image of $\phi(F[\epsilon])$. Since $F[\epsilon]$ is local, Grothendieck-Hilbert's Theorem 90 yields $H^1(F[\epsilon],\G_m)=0$, so that  $f$ lifts to $$f' \in \GL(\Gamma^n(V))(F[\epsilon]).$$ Since $\psi(f)=\Id$, there exists a morphism of $F$-schemes \[\lambda: (\A(V) -\{ 0 \}) \to \R_{F[\epsilon]/F}(\Gm) \simeq \G_m \times_F \A^1,\] \[ v \mapsto  \lambda_1(v)+ \lambda_2(v) \epsilon,\]such that $$ f'( [v]_n)= \lambda(v) [v]_n, $$ identically on points. Let us check that $\lambda$ is constant. Since $\dim(V) \geq 2$, and the source of $\lambda$ is normal, it extends to a morphism of $F$-varieties \[\Lambda: \A(V)  \to \G_m \times_F \A^1 \subset \A^2.\]   Denote by $\delta$ the degree of $\Lambda$, as a polynomial map. From the equality $$ f'( [v]_n)= \lambda(v) [v]_n, $$ also valid on functors of points, we get $n=\delta+n$, whence $\delta=0$ and $\Lambda$ is constant. Rescaling $f'$, we can thus assume $\Lambda=1$. Since $F$ is infinite, pure symbols $[v]_n$ span the $F$-vector space $\Gamma^n(V)$, so that $f'=\Id$. Hence, $f=\Id$, as desired. 
\end{dem}

\subsection{Characteristic-free polarity.}

\begin{defi}\label{DefiGeoDivPow}
     Let $W$ be an $F$-vector space, let $  X \subset \P(W)$ be a closed $F$-subscheme, defined by a sheaf of ideals $$0 \to \mathcal I_X \to \mathcal O_{\P(W)} \to \mathcal O_X \to 0.$$ Define $$ E_{X,m}:=  H^0(\P(W), \mathcal I_X(m))\subset  H^0(\P(W), \mathcal O(m))=\Sym^m(W ^\vee).$$ For all sufficiently large $m$, it generates $\mathcal I_X(m) $.\\ For brevity, denote $E_{X,m}$ by $E_X$. Dualizing the exact sequence $F$-vector spaces $$ 0 \to E_X \to \Sym^m(W^\vee) \to \ast \to 0,$$ one gets an exact sequence $$ 0 \to L_X \to \Gamma^m(W) \to E_X^\vee \to 0.$$ 
\end{defi}
\begin{lem}\label{LemGeoDivPow}
  Keep notation of Definition \ref{DefiGeoDivPow}. For $m$ large enough, there is a natural isomorphism of $F$-schemes $$ P_X: X \to \P(L_X) \bigcap \Ver_m(\P(W))$$ $$ z \mapsto [z]_m,$$  where $\bigcap$ denotes scheme-theoretic intersection in $\P(\Gamma^m(W))$.
\end{lem}

\begin{dem}
  Let $R$ be an $F$-algebra. Since $m>>0$, the set $X(R) \subset \P(W)(R)$ consists of those lines, on which all $m$-linear forms in $E_X$ vanish. Using duality between $\Sym^m(W^\vee)$ and $\Gamma^m(W)$, this translates as $$ X(R)=\{ (w) \in \P(W)(R),\;  \phi([w]_m)=0,\;  \forall \phi \in E_W\}.$$ Via the closed immersion $\Ver_m$, the right side of the equality coincides with  $ \P(L_X)(R) \bigcap\Ver_m(\P(W))(R) $. This holds for any $R$, whence the desired iso $P_X$.
\end{dem}

\subsection{Morphisms of varieties induced by multiplication of $\Gamma(V)$.}

The following notion is especially important  if $\carac(F)=p$.
\begin{defi}($F$-disjointness).\\
Let $ a_1,a_2, \ldots a_d$ be  positive integers.\\ Say that $a_1,a_2,\ldots,a_d$ are $F$-disjoint if the following holds. \begin{enumerate}
    
    \item{ For  $i=1,\ldots, d-1$, one has $a_i+a_{i-1}+\ldots+ a_1 < a_{i+1}$.} \item{If $\carac(F)=p$, for  $i=1,\ldots, d-1$ one has $a_i+a_{i-1}+\ldots+ a_1 < p^{v_p(a_{i+1})} $.}
\end{enumerate}
\end{defi}
\begin{ex}
    Assume  $\carac(F)=p$, and  $a_i=p^{r_i}$, with $0 \leq r_1 < r_2 < \ldots < r_d$. Then, $a_1,a_2,\ldots,a_d$ are $F$-disjoint.
\end{ex}
\begin{rem}
    If  $\carac(F)=p$, then $(2)$ implies $(1)$ in Definition above. Thinking in  base $p$,   $(2)$ is equivalent to the following.  The position of the least non-zero digit of $a_{i+1}$, is stricly bigger than that of the  greatest non-zero digit of $a_i+a_{i-1}+\ldots+ a_1$.
\end{rem}
Recall a well-known fact.
\begin{lem}\label{LemDisjoint}
	
 Let $ a_1, \ldots, a_d$ be  nonnegative integers.\\ The  $p$-adic valuation of the multinomial coefficient $ { a_1+ \ldots + a_d  \choose {a_1, \ldots, a_d}}$ is the number of carryovers, when computing the  sum $a_1+  \ldots+ a_d$ in base $p$. \\ Hence, if $\carac(F)=p$ and $a_1,a_2,\ldots,a_d$ are $F$-disjoint, $ { a_1+ \ldots + a_d  \choose {a_1, \ldots, a_d}}$ is prime-to-$p$. 
 \end{lem}

\begin{dem}
One can apply induction on  $r \geq 2$, using the formula $$ { a_1+ \ldots + a_d \choose {a_1, \ldots, a_d}} = { a_1+ \ldots + a_d \choose {a_1+a_2, a_3, \ldots, a_d}}  {a_1 + a_2 \choose a_1, a_2 }.$$ The claim to prove when $r=2$ is a classical fact, which is also a nice elementary exercise left to the reader. The second assertion  readily follows.
\end{dem}

\begin{lem}\label{LemDivProd}
    Let $V$ be an $F$-vector space.   The following is true. \begin{enumerate}
        \item{Let $a,b$ be $F$-disjoint  integers, and let  $y \in V-\{0 \}$. The  multiplication    $$M_y: \Gamma^a(V)   \to \Gamma^{a+b}(V) $$ $$x \mapsto x [y]_b  $$ is an $F$-linear injection.}  \item{The  formula $$\mu: \P(\Gamma^a(V)) \times_F \P(V)   \to \P(\Gamma^{a+b}(V)) $$ $$(x,y) \mapsto x [y]_b  $$ defines a morphism of $F$-varieties, injective on $\overline F$-points. }\item{ Let $ a_1,a_2, \ldots a_d$ be $F$-disjoint  integers. Then, the  formula
    $$ \tau: \P(V) \times_F  \ldots \times_F \P(V)  \to \P(\Gamma^{a_1+\ldots+a_{d}}(V))  $$ $$(x_1, x_2, \ldots, x_d) \mapsto [x_1]_{a_1} [x_2]_{a_2} \ldots[x_d]_{a_d}$$ defines a morphism of $F$-varieties, injective on $\overline F$-points.}
    \end{enumerate}
    
\end{lem}

\begin{dem}
Can assume  $F=\overline F$.  Let us prove item (1). Pick a basis $(e_1, \ldots, e_n)$ of $V$, with $e_n=y$. Work in the standard basis $[e_1]_{a_1} \ldots [e_n]_{a_n}$ of $\Gamma^a (V)$, indexed by partitions $a=a_1+ \ldots+ a_n$. Similarly, work in the standard basis $[e_1]_{c_1} \ldots [e_n]_{c_n}$ of $\Gamma^{a+b} (V)$, indexed by partitions $a+b=c_1+ \ldots+ c_n$. Let us compute: $$ ([e_1]_{a_1} \ldots [e_n]_{a_n}) [y]_b={a_n +b \choose  b } [e_1]_{a_1} \ldots [e_{n-1}]_{a_{n-1}}  [e_n]_{a_n+b}. $$  If $\carac(F)=0$, it readily follows that $M_y$ is injective. Assume  $\carac(F)=p$. \\ Since $a_n \leq a< p^{v_p(b)}$, computing $a_n+b$ in base $p$ occurs without carryovers. Thanks to Lemma \ref{LemDisjoint},  ${a_n +b \choose  b } \in F$ is non-zero. Consequently, $M_y$ is still injective. \\
Let us prove that $(2)$ implies $(3)$. If $d \geq 3$, $\tau$ factors as the composite of $$\P(V)   \times_F \ldots \times_F \P(V)\times_F \P(V)  \to \P(\Gamma^{a_1+\ldots+a_{d-1}}(V)) \times_F \P(V) $$  $$(x_1, \ldots, x_d) \mapsto ([x_1]_{a_1}\ldots  [x_{d-1}]_{a_{d-1}},x_d)$$ and  $$\P(\Gamma^{a_1+\ldots+a_{d-1}}(V)) \times_F \P(V)  \xrightarrow{\mu} \P(\Gamma^{a_1+\ldots+a_d}(V))$$  $$(x,y) \mapsto x [y]_{a_d} .$$ By induction, item (3) thus  indeed follows from (2).\\
It remains to prove (2).  That $\mu$ is well-defined, follows from (1) (injectivity of $M_y$, for $y \neq 0$). Let us check injectivity of  $\mu$ on $F$-points. Let $y,y' \in V- \{0 \}$ and $x,x' \in \Gamma^a(V)- \{0 \}$, be such that $\mu(x,y)=\mu(x',y')$. Rescaling, one can assume $$x [y]_b= x' [y']_b \in \Gamma^{a+b}(V).$$ Suppose that $(y) \neq (y') \in \P(V)(F)$. 
 Pick a basis $(e_1=y, e_2=y',e_3,\ldots, e_n)$ of $V$. Working in the monomial basis $([e_1]_{a_1} \ldots [e_n]_{a_n})$ of $\Gamma^a (V)$, and using $a <b$, one sees that  $\Im(M_y) \bigcap \Im(M_{y'})= \{0\}$, contradicting $M_y(x) =M_{y'}(x')\neq 0$. Hence $y$ and $y'$ are collinear. Rescaling them, one can assume $y=y'$. Using item (1), one  concludes that $x=x'$, which finishes the proof.

\end{dem}

\begin{rem}
    In general,  morphisms in items (2) and (3) above, are not injective on tangent spaces: they are not closed immersions.
\end{rem}
\begin{rem}
    If $\carac(F)=0$, one can then replace $\Gamma^a(V)$ by $\Sym^a(V)$, and accordingly replace symbols $[x]_a$ by $\frac {x^a} {a!}$. Using  that the polynomial $F$-algebra $\Sym(V)$ is a UFD, the proof of Lemma \ref{LemDivProd} is then   easier.
\end{rem}
\section{Concrete Tannakian construction.}

\subsection{Explicit action of $\PGL_d$ with trivial stabilisers.}

\begin{lem}\label{LemFreeDiv}

Let $d \geq 3$,  and let $V$ be a $d$-dimensional $F$-vector space. \\Let $1,a_1,a_2, \ldots, a_{d+1}$ be $F$-disjoint integers. Set $r:=a_1+ \ldots +a_{d+1}$. \\
Choose a basis $(e_1,\ldots, e_d)$ of $V$, and define $$ r:=a_1+\ldots+a_{d+1},$$ $$e_{d+1}:=e_1+  e_2+ \ldots+ e_d, $$ $$ x:=[e_1]_{a_1}[e_2]_{a_2} \ldots [e_d]_{a_d} [e_{d+1}]_{a_{d+1}} \in \Gamma^r(V).$$ Then $x \neq 0$, and $$\Stab_{\PGL(V)}((x) )=\{1\},$$ for the natural action of $\PGL(V)$ on $\P(\Gamma^r(V))$.
\end{lem}

\begin{dem}

That $x \neq 0$ follows from item (3) of Lemma \ref{LemDivProd}. There is an extension of $F$-groups  $$1 \to \mu_r \to \Stab_{\GL(V)}(x) \to \Stab_{\PGL(V)}((x))  \to 1.$$ Note that its kernel is not étale if $\carac(F)=p$. Let us show the triviality of $\Stab_{\PGL(V)}((x))$. To do so, one may assume $F=\overline F$.\\ Pick $f \in \GL(V)(F)$. Assume that $f(x)=x$. Using the injectivity statement of  Lemma \ref{LemDivProd} (3), one sees that     $f$ fixes each $e_i$ up to scalars, implying that $f$ is homothetic.  In other words, the group  $\Stab_{\PGL(V)}((x))(F)$ is trivial. If $\carac(F)=0$, this finishes the proof. If $\carac(F)=p$,  it remains to prove triviality of $\Lie(\Stab_{\PGL(V)}(x))$.  \\In  computations that will follow, one typically uses the formula $$ [x]_a [x]_b= {{a+b} \choose a} [x]_{a+b},$$ for various $a,b \in \N$, and one checks whether or not ${{a+b} \choose a} $ is divisible by $p$.\\ Observe that $p$ divides $a_1, \ldots, a_{d+1}$, because $1,a_1,a_2, \ldots, a_{d+1}$ are $F$-disjoint. \\Pick $u \in \End(V)$, with matrix $(u_{i,j})_{1 \leq i,j \leq d}$ in the basis $(e_1,\ldots, e_d)$. \\Set $u_i:=u(e_i)$. If $i \leq d$, then $u_i=\sum_{j=1}^d u_{j,i}e_j$, and $ u_{d+1}= u_1+ \ldots +u_d.$ Set $$f:=(\Id+\epsilon u)  \in \GL(V)(F[\epsilon]).$$ Assume there exists $c \in F$, such that  $$ f(x)=(1+c \epsilon)x.$$ To conclude the proof, one  needs to  show  $u\in F \Id$. \\Note that $f(x)$ reads as  $$ x=[e_1+\epsilon u_1]_{a_1}[e_2+\epsilon u_2]_{a_2} \ldots [e_{d+1}+\epsilon u_{d+1}]_{a_{d+1}} .$$ Expanding divided powers, and comparing coefficients of $\epsilon$, one gets $$(E):c [e_1]_{a_1} \ldots[e_{d+1}]_{a_{d+1}}=\sum_{i=1} ^{d+1}[u_i]_1 [e_1]_{a_1} [e_2]_{a_2}\ldots [e_i]_{a_i -1} \ldots [e_{d+1}]_{a_{d+1}}  .$$  In the monomial basis of $\Gamma^r(V)$ furnished by $e_1,\ldots, e_d$, consider the coefficient  of $$M:=[e_1]_{a_1-1} [e_2]_{a_2+1} [e_3]_{a_3+a_{d+1}}[e_4]_{a_4} [e_5]_{a_5} \ldots [e_d]_{a_d},$$ of both sides of this equality. It vanishes on the left side. In  the sum on the right side,  only $i=1$ can contribute. Let us  expand the  corresponding term, reading as  $$ [u_{1,1} e_1+\ldots + u_{d,1} e_d]_1 [e_1]_{a_1-1} [e_2]_{a_2}[e_3]_{a_3} \ldots [e_d]_{a_d}  [ e_1+  e_2+ e_3+ \ldots+e_d]_{a_{d+1}}.$$ In the decomposition $$  [ e_1+  e_2+ e_3+ \ldots+e_d]_{a_{d+1}}=\sum_{b_1+ \ldots+ b_d=a_{d+1} } [e_1]_{b_1} [e_2]_{b_2} [e_3]_{b_3} \ldots [e_d]_{b_d}, $$ the only two partitions that may contribute to a non-zero multiple of $M$, are $$(b_1,b_2,b_3,b_4, \ldots, b_d)=(0,0,a_{d+1},0, \ldots 0) $$ and $$(b_1,b_2,b_3,b_4, \ldots, b_d)=(0,1,a_{d+1}-1,0, \ldots 0).$$ These terms are given, respectively, by $$ u_{2,1} [e_2]_1 [e_1]_{a_1-1} [e_2]_{a_2}[e_3]_{a_3} \ldots [e_d]_{a_d}  [  e_3]_{a_{d+1}}={{a_3+ a_{d+1}} \choose {a_3}} (a_2+1)u_{2,1} M,$$ and $$ u_{3,1} [e_3]_1 [e_1]_{a_1-1} [e_2]_{a_2}[e_3]_{a_3} \ldots [e_d]_{a_d}  [e_2]_1 [  e_3]_{a_{d+1}-1}={{a_3+ a_{d+1}} \choose {1,a_3,a_{d+1}-1}} (a_2+1)u_{3,1}M.$$ Gathering the information above, one gets $$0={{a_3+ a_{d+1}} \choose {a_3}} (a_2+1)u_{2,1}+{{a_3+ a_{d+1}} \choose {1,a_3,a_{d+1}-1}} (a_2+1)u_{3,1} .$$ Since $1,a_1,a_2, \ldots, a_{d+1}$ are $F$-disjoint, Lemma \ref{LemDisjoint} asserts  that $p$ does not divide  ${{a_3+ a_{d+1}} \choose {a_3}} (a_2+1)$, but divides ${{a_3+ a_{d+1}} \choose {1,a_3,a_{d+1}-1}}$ (for the latter fact, observe that adding $1$ and $(a_{d+1}-1)$ in base $p$, occurs with carryovers). Thus $u_{2,1}=0$. One can reproduce this argument, with any triple of distinct indices $\in \{1, \ldots,d\}$, in place of $(1,2,3)$. One thus  gets  $u_{i,j}=0$ for all $i \neq j$. Thus, $u_i=\alpha_i e_i$,  $i=1, \ldots, d$. In $(E)$, put the term of index $i=d+1$  on the other side of the equation. This gives $$(E'): c [e_1]_{a_1} \ldots[e_{d+1}]_{a_{d+1}} -[\alpha_1 e_1+ \ldots\alpha_d e_d ]_1 [e_1]_{a_1} \ldots  [e_d]_{a_d }  [e_{d+1}]_{a_{d+1}-1} $$ $$=\sum_{i=1} ^d a_i\alpha_i [e_1]_{a_1} [e_1]_{a_2}\ldots [e_i]_{a_i} \ldots [e_{d+1}]_{a_{d+1}}=0.$$ To finish, let us work in the  following basis  of $V$: $$(f_1:=-e_2, f_2:=-e_3, \ldots, f_{d-1}:=-e_d,f_d:= e_{d+1}),$$ and in the induced monomial basis of $\Gamma^r(V)$. Note that $$ e_1=  f_1+ \ldots +  f_d.$$Equality $(E')$ gives $$c [e_1]_{a_1} [f_1]_{a_2} \ldots[f_d]_{a_{d+1}} =[\alpha_1 e_1+ \ldots\alpha_d e_d ]_1 [e_1]_{a_1} [f_1]_{a_2} \ldots  [f_{d-1}]_{a_d }  [f_d]_{a_{d+1}-1}. $$ Express both sides in the  monomial basis,  and consider the coefficient of $$[f_1]_{a_2+1} [f_2]_{a_3} \ldots[f_d]_{a_{d+1}}. $$ On the left side, it is $c(a_2+1)=c$. Since $$[f_d]_1 [f_d]_{a_{d+1}-1}=a_{d+1}[f_d]_{a_{d+1}}=0,$$ it is $0$ on the right side, so that  $c=0$. Thus, $$(E''): [\alpha_1e_1+ \ldots\alpha_d e_d ]_1 [e_1]_{a_1} [f_1]_{a_2} \ldots  [f_{d-1}]_{a_d }  [f_d]_{a_{d+1}-1}=0. $$ Let $\beta_i \in F$ be such that $$\alpha_1e_1+ \ldots\alpha_d e_d= \beta_1 f_1+ \ldots +\beta_d f_d.$$   Computing the coefficient of $[f_1]_{1+a_1+a_2} [f_2]_{a_3} [f_3]_{a_4} \ldots [f_{d-1}]_{a_d}[f_d]_{a_{d+1}-1}$ in $(E'')$, one gets $$ \beta_1  { 1+a_1 + a_2  \choose {1,a_1,  a_2}}=0 \in F,$$ where the multinomial coefficient  is prime-to-$p$ by Lemma \ref{LemDisjoint}. Thus $\beta_1=0$.\\ In the same fashion,  $\beta_2= \ldots= \beta_{d-1}=0$. In other words: all $\alpha_i$ are equal to $\beta_d$, so that $u=\beta_d \Id$, as was to be shown.\\

\end{dem}

\subsection{Linear algebraic groups as stabilisers of symbols.}

\begin{prop}\label{PropZ}

Let  $G$ be a linear algebraic group over $F$. There exists a  representation $G \hookrightarrow \GL(W)$, such that the composite $G \hookrightarrow \GL(W) \to \PGL(W)$ is  faithful, together with the following data. \begin{enumerate}
    \item{A closed subscheme $Z \subset \P(W)$, such that $G=\Stab_{\PGL(W)}(Z)$.} \item{A $G$-fixed rational point $(w_0) \in \P(W)(F)-Z(F)$.}
\end{enumerate}

\end{prop}
\begin{dem}
    Pick  $n \geq 2$ and a faithful  representation $G \hookrightarrow \GL_{n-1}$. Note that the natural composite  $$G \hookrightarrow  \GL_{n-1}\subset \GL_n\to  \PGL_n$$  is still an embedding. This way, one gets a faithful representation $G \hookrightarrow \PGL_n$, such that the action of $G$ on $\P^{n-1}$ has $(e_n)$ as an $F$-rational fixed point. Define $$d:=2n,  \; V_1=V_2:=F^n, \; V:=V_1 \bigoplus V_2.$$ Pick $n$ large enough, so that $d-1> \dim(G)$. Consider the diagonal composite  $$ G \hookrightarrow \GL(V_1) \stackrel {x \mapsto (x,x)} \hookrightarrow \GL(V),$$ inducing $$ G \hookrightarrow \PGL(V_1) \stackrel {x \mapsto (x,x)} \hookrightarrow \PGL(V).$$Consider the canonical basis $(e_1, \ldots, e_n, e_{n+1}, \ldots, e_{2n})$ of $V$, obtained by putting together two copies of the canonical basis of $V_1$.  Let $(a_1,a_2, \ldots, a_{d+1})$, $r=a_1+\ldots +a_{d+1}$ and $x \in \Gamma^r(V)$ be as in the premises of Lemma \ref{LemFreeDiv}. Setting $W:=\Gamma^r(V)$, this Lemma states that  $\Stab_{\PGL(V)}((x) \in \P(W)(F) )$ is trivial.  \\  Since $\PGL(V)$ is a smooth $F$-group, acting on the smooth $F$-variety $\P(W)$, it is known that the  $\PGL(V)$-orbit $$ O(\PGL(V),x):=\PGL(V).x$$ is a locally closed subscheme of $\P(W)$. Indeed, its closure $$\overline{ O(\PGL(V),x)} \subset \P(W)$$ equipped with its reduced induced scheme structure, is a $\PGL(V)$-stable closed subscheme, and  the orbit $O(\PGL(V),x) \subset \overline{ O(\PGL(V),x)}$ is open in its closure. Denote by $Z\subset \overline{ O(\PGL(V),x)}$ its  complement, considered with its reduced induced structure. Set $U:= \P(W)-Z$. One has $G$-equivariant embeddings $$G \stackrel {\mathrm {closed}} \hookrightarrow \PGL(V)  \stackrel {\xrightarrow{g \to g.x}} \sim O(\PGL(V),x) \stackrel {\mathrm {closed}} \hookrightarrow U   \stackrel {\mathrm {open}}\hookrightarrow \P(W)$$ Thus, the $G$-orbit map $$ \alpha: G  \xrightarrow{g \to g.x}\P(W)$$ is a locally closed immersion (even though $G$ may not be smooth).  Consider its scheme-theoretic image $X \stackrel {\mathrm {closed}} \hookrightarrow \P(W)$. Let us check that $X$ is $G$-stable. Let $R$ be an $F$-algebra. Since formation of scheme-theoretic image of a quasi-compact morphism commutes to flat base-change (\cite{SP}, 100.38.5), the $R$-scheme $X_R:=X \times_{\Spec(F)} \Spec(R)$ is the scheme-theoretic image of the $G_R$-orbit map $\alpha_R$. Pick $\gamma \in G(R)$. Then $\gamma. X_R$ is the scheme-theoretic image of the $R$-morphism $$ (\gamma.\alpha_R): G_R  \xrightarrow{g \to (\gamma g).x}\P(W)_R.$$  This morphism factors as $$ G_R   \xrightarrow{g \to \gamma g } G_R  \xrightarrow{\alpha_R}\P(W)_R. $$ Since $g \mapsto \gamma g$ is an isomorphism, we conclude that $(\gamma.\alpha_R)$ and $\alpha_R$ share the same scheme-theoretic image. Equivalently, $\gamma. X_R=X_R$, proving that $X$ is $G$-invariant.
Let us check that inclusion of linear algebraic $F$-groups $$G \subset \Stab_{\PGL(V)} (X)$$ is  an equality.  Set $$Y:=X - O(G,x) \subset \P(W).$$ It is a closed subset of $X$. Consider it as a closed subscheme of $X$, using the reduced induced structure (as such, it may not be $G$-invariant).  Let $R$ be a finite local $\overline F$-algebra, and let $\phi \in \Stab_{\PGL(V)} (X)(R)$. Arguing by contradiction,  suppose that $\phi.x \notin O(G,x)(R).$ Denote by $\phi_0 \in \Stab_{\PGL(V)} (X)(\overline F)$ the special fiber of $\phi$. Since $R$ is local and $ O(G,x) \subset X$ is open, one has $\phi_0.x \notin O(G,x)(\overline F).$ \\The  monomorphism of $F$-schemes$$\beta: G \to X$$ $$ g \mapsto (g \phi_0).x $$ would then, set-theoretically, take values   in $Y$. Indeed, suppose that there exists  $g \in G(\overline F)$, such that  $(g \phi_0).x \in O(G,x)(\overline F)$. Because $O(G,x)$ is a principal homogeneous space of $G$, there exists $\gamma \in G(\overline F),$ such that $(g \phi_0).x=\gamma.x$. Then $$\gamma^{-1} g \phi_0 \in \Stab_{\PGL(V)}(x)(\overline F)=\{1\},$$  implying $\phi_0 \in G(\overline F)$, hence  $\phi_0.x \in O(G,x)(\overline F),$ contradicting $\phi_0.x \notin O(G,x)(\overline F).$ Thus, set-theoretically, the monomorphism $\beta$  takes values in $Y$. This is impossible because $\dim(G) > \dim(Y)$, as Noetherian topological spaces. One concludes that $\phi.x \in O(x)(R)$, implying  $\phi \in G(R)$. This proves  $G =\Stab_{\PGL(V)} (X)$.\\  
It remains to prove that $X$ does not intersect $\Ver_r(\P(V) \subset \P(W).$    Recall that $V=V_1 \bigoplus V_2$, $V_1=V_2$, and that the $G$-action on $\P(V)$ occurs via the  composite $$ G \hookrightarrow \PGL(V_1) \xrightarrow{diag} \PGL(V).$$ By choice of the basis $(e_1, \ldots, e_d)$, the orbit $O(\PGL(V_1), x)$  is thus contained in the image of the composite $F$-morphism  $$ \pi: \P(V_1)^n  \times_F  \P(V_2)^n  \times_F \P(V)  \hookrightarrow \P(V)^n  \times_F  \P(V)^n  \times_F \P(V) \xrightarrow{\tau} \P(\Gamma^{a_1+\ldots+a_{d+1}}(V) ) $$ $$(x_1, x_2, \ldots, x_{d+1}) \mapsto [x_1]_{a_1} [x_2]_{a_2} \ldots[x_{d+1}]_{a_{d+1}}.$$  Here the first arrow is obtained by taking products of the natural closed immersions $\P(V_i) \hookrightarrow \P(V)$, $i=1,2$. The arrow $\tau$ is that of item (3) of Lemma \ref{LemDivProd}. \\Observe  that  $\Ver_r$ is the composite $F$-morphism $$ \P(V) \xrightarrow{diag} \P(V)^{d+1} \xrightarrow{\tau} \P(\Gamma^{a_1+\ldots+a_{d+1}}(V) ),$$ $$ x \mapsto [x]_{a_1} [x]_{a_2} \ldots[x]_{a_{d+1}}. $$  By item (3) of Lemma \ref{LemDivProd}, $\tau$ is injective on $\overline F$-points. Since $\P(V_1)$ and $\P(V_2)$ intersect trivially as linear subspaces of $\P(V)$, it follows  that  $$ \Im(\pi)(\overline F)\bigcap    \Ver_r(\P(V))(\overline F) =\varnothing.$$\\ Since the  source of $\pi$ is proper, one gets $$\overline{O(\PGL(V_1), x)}(\overline F) \subset \Im(\pi)(\overline F),$$ so that  $$ \overline{O(\PGL(V_1), x)} \bigcap \Ver_r(\P(V)) =\varnothing,$$ and a fortiori $$X \bigcap \Ver_r(\P(V))= \varnothing.$$
Next, consider  the closed subscheme $$Z:=X \coprod \Ver_r(\P(V)) \subset \P(W),$$  which is indeed a disjoint union.   Let us check that the natural embedding $$ G \hookrightarrow \Stab_{\PGL(W)}(Z)$$ is an iso. To do so, let $R$ be a finite local  $\overline F$-algebra, and let  $$\phi \in \Stab_{\PGL(W)} (Z)(R).$$  Then $\phi(\Ver_r(\P(V)) _R) \subset Z_R$ is an irreducible smooth clopen $R$-subscheme, of dimension $d-1 > \dim(G)=\dim(X)$. It thus intersects $X_R \subset Z_R$  trivially. In other words:  $\phi(\Ver_r(\P(V))_R) =\Ver_r(\P(V))_R$, and consequently $\phi(X_R) =X_R$.\\By Proposition \ref{AutVeronese}, $\phi$ belongs to $\PGL(V)(R)$. Since $G= \Stab_{\PGL(V)}(X)$, one then gets $\phi \in G(R)$, as was to be shown. Item (1) is proved.\\ For (2), recalling that $e_n \in H^0(G,V_1)$ and $e_{2n} \in H^0(G,V_2)$, one may take $$(w_0):= [e_n]_{a_1} [e_{2n}]_{r-a_1}=\tau(e_n, e_{2n}, e_{2n}, \ldots, e_{2n}) \in H^0(G,\P(W)(F)).$$ The injectivity of $\tau$ (already used above to prove $X \bigcap \Ver_r(\P(V))= \varnothing$) then guarantees that $(w_0) \notin \Ver_r(\P(V))(\overline F)$ and $(w_0) \notin \Im(\pi)(\overline F)$. Hence  $(w_0) \notin Z(F)$.

\end{dem}

\begin{prop}\label{DividedTannaka}
Let  $G$ be a linear algebraic group over $F$. There exists an $F$-vector space $W$,   integers $n,l \geq 1$, and a  linear subspace  \[ L \in \Gr(l, \Gamma^n(W))(F),\] such that \[ G \stackrel \sim \to \Stab_{\PGL( W)}(L) ,\] as group schemes over $F$. \\ Moreover, one can take $L$ such that the closed subvarieties $\P(L)\stackrel {\mathrm{lin}} \hookrightarrow \P( \Gamma^n(W))$ and $\P(W) \stackrel {\mathrm {Ver_n}} \hookrightarrow \P(\Gamma^n(W))$ do not intersect.

\end{prop}
\begin{dem}
Pick a $G$-representation $W$, a closed subvariety $Z \subset \P(W)$ and $(w_0) \in \P(W)(F)$  as in Proposition \ref{PropZ}. For $m \in \N $,  consider the $F$-subspaces $E_Z \in \Gr(l,\Sym^m(W^\vee))(F)$ and $ L_Z\in \Gr(l,\Gamma^m(W))(F),$ introduced in Definition \ref{DefiGeoDivPow}. These are $G$-stable. Fix $m$ large enough, so that  $m \neq -1 \in F$, and  $$ G= \Stab_{\PGL(W)}(E_Z  )$$  (see Lemma \ref{LemGeoDivPow}). Considering the  exact sequences of  Definition \ref{DefiGeoDivPow}, one sees that $$ \Stab_{\PGL(W)}(E_Z  )=\Stab_{\PGL(W)}(E_Z^\vee  )=\Stab_{\PGL(W)}(L_Z  ),$$ so that $$ G= \Stab_{\PGL(W)}(L_Z  ).$$  If $\carac(F)=0$ (resp. $\carac(F)=p$), set $q:=m+1$ (resp. $q:=p^s>m$,  the smallest $p$-th power greater than $m$). Set $n:=m+q$. Consider the $F$-linear map $$M_{w_0}: \Gamma^m(W) \to \Gamma^n(W)$$  $$x \mapsto x[w_0]_q. $$ It is injective by Lemma \ref{LemDivProd}. Set  $$ L:= M_{w_0}(L_Z)  \subset  \Gamma^n(W).$$ Since $(w_0) \in H^0(G,\P(W))$, there is a natural inclusion of $F$-groups $$  \Stab_{\PGL(W)}(L_Z)=G\subset  \Stab_{\PGL(W)}(L) .$$ Let us show it is an equality. To do so, one may assume $F=\overline F$. Let $A$ be a  finite local $F$-algebra with maximal ideal $\mathcal M$, and let $g \in  \Stab_{\PGL(W)}(L)(A)$. We need to show that $g \in  \Stab_{\PGL(W)}(L_Z)(A)$.  Let us first show that $g$ fixes $(w_0) \in \P(W)(F) \subset \P(W)(A)$. Assume that $A=F$. If $w_0$ and $g(w_0)$ are not $F$-collinear, complete them into an $F$-basis $(w_0,w_1=g(w_0),w_2, \ldots, w_d)$ of $W$. By assumption, the two subspaces $$ g(L)=[w_1]_q. g(L_Z) \subset \Gamma^{n}(W)$$ and $$ L=[w_0]_q. L_Z \subset \Gamma^{n}(W)$$ are equal. Work in the natural basis of $\Gamma^{n}(W)$ induced by $(w_0,w_1,w_2, \ldots, w_d)$. Then, elements of $ L$ are linear combinations of symbols of the shape $$(A): [w_0]_{a_0}\ldots [w_d]_{a_d},$$ with $a_0 \geq q$ and $a_0+ \ldots +a_d=n <2q$. Similarly,   elements of $ g(L)$ are linear combinations of symbols of the shape $$(B): [w_0]_{b_0} [w_1]_{b_1}\ldots [w_d]_{b_d},$$ with $b_1 \geq q$ and $b_0+ \ldots +b_d=n<2q$. But no symbol  is  of both  shapes (A) and (B)- a contradiction. Hence $g$ fixes $(w_0)$. If $\carac(F)=0$, this is enough to conclude. It remains to treat the case $\carac(F)=p$, $q=p^s >m$ and $A$  arbitrary. Denote by $\overline g \in  \Stab_{\PGL(W)}(L)(F)$ the residue of $g$. Pick $z \neq 0 \in W$, such that $(z) \in Z(F)$, so that $[z]_m \in L_Z$ and $(z) \neq (w_0)$. Define $$w_1:=\overline g (z) \in W.$$ By the case $A=F$  dealt with before, $$(\overline g(w_0))=(w_0) \in \P(W)(F),$$ so that $(w_0) \neq (w_1)$. Complete $w_0,w_1$ into an $F$-basis $(w_0,w_1,w_2, \ldots, w_d)$ of $W$. Rescaling $g$ by an element of $A^\times$, one can assume $$ g(w_0)=w_0 + \epsilon_1 w_1+ \ldots+ \epsilon_d w_d,$$ where $\epsilon_i \in \mathcal M$. Assume first, that $\epsilon_i\mathcal M=0$ for $i=1, \ldots, d$.\\  There exists $\eta \in \Gamma^m(W)\otimes_F \mathcal M $ such that  $$  [g(z)]_m=[w_1]_m + \eta \in \Gamma^m(W)\otimes_F A,$$ and  a little computation   gives $$ [g(w_0)]_q=[w_0]_q+ \epsilon_1[w_0]_{q-1}[w_1]_1+ \ldots+ \epsilon_d[w_0]_{q-1}[w_d]_1\in \Gamma^q(W)\otimes_F A. $$ Developping the product, rearranging terms, one gets $$[g(w_0)]_q [g(z)]_m= (m+1)\epsilon_1[w_0]_{q-1}[w_1]_{m+1}+\sum_{i=2}^d\epsilon_i[w_0]_{q-1}[w_1]_m [w_i]_1 + [w_0]_q E,$$ for some $E \in \Gamma^m(W)\otimes_F A$. Recall that all elements of $L$ are linear combinations of symbols of shape (A) above. Observe that symbols $[w_0]_{q-1}[w_1]_{m+1}$ and $[w_0]_{q-1}[w_1]_m [w_i]_1$ are not of shape (A). Since $L=g(L)$, it must be the case that $\epsilon_i=0$ for $i=2, \ldots, d$ . Since $m \neq -1 \in F$, one also has $\epsilon_1=0$.\\
It remains to remove the assumption $\epsilon_i\mathcal M=0$. This is a straightforward induction on $k\geq 1$, such that $\mathcal M^k=0$. If $k =1 $ there is nothing to do. Assume that  $\mathcal M^{k+1}=0$. By induction applied to $A/\mathcal M^k$, one gets $\epsilon_i \in \mathcal M^k$, so that  $\epsilon_i\mathcal M=0$ and the above applies, yielding $\epsilon_i=0$.\\
We have proved $g((w_0))=(w_0) \in \P(W)(A)$. \\By  item (1) of Lemma \ref{LemDivProd}, the $A$-linear map $$M_{w_0}: \Gamma^m(W)\otimes_F A \to \Gamma^{m+q}(W)\otimes_F A $$  $$x \mapsto x[w_0]_q $$ is injective. Since  $$M_{w_0}(L_Z)=M_{g(w_0)}(g(L_Z))=M_{w_0}(g(L_Z)),$$ it is then straightforward to see that $g(L_Z)=L_Z$. We have proved $$ G= \Stab_{\PGL(W)}(L_Z)=\Stab_{\PGL(W)}(L) .$$
To conclude, it remains to prove that $\P(L)$ and $\Ver_n(\P(W))$ intersect trivially. \\By item (2) of Lemma \ref{LemDivProd}, the morphism $$\mu:  \P(\Gamma^m(V)) \times_F \P(V) \to \P(\Gamma^n(V)) , $$ $$(x,y) \mapsto x[y]_q  $$ is injective on $\overline F$-points. Introduce the graph of $\Ver_m$, $$\Delta : \P(V) \to \P(\Gamma^m(V)) \times_F \P(V),$$ $$ v \mapsto ([v]_m,v).$$ Observe that  $$\P(L)=\mu (\P(L_Z) \times \{w_0 \})$$ and $$\Ver_n(\P(W))=\mu (\Delta(\P(V))).$$  Because $(w_0) \notin Z(\overline F)$, one has $[w_0]_m \notin L_Z$, so that  $$(\P(L_Z) \times \{w_0 \}) \bigcap \Delta(\P(V)) =\varnothing. $$ Thus, $\P(L) \bigcap  \Ver_n(\P(W))= \varnothing$, as desired.

\end{dem}

\begin{qu}\label{Ql}
    In Proposition \ref{DividedTannaka},  can one take $l=1$ ? 

\end{qu}
 
\begin{rem}
    We suspect that the answer to Question \ref{Ql} is yes, and sketch an optimistic strategy to investigate it.\\ Let $W,l,n$ and $L$ be furnished by Proposition \ref{DividedTannaka}. By inspection of its proof, one may assume that $n$  is odd.  Then, there is a well-defined $F$-linear map $$\Psi:  \Lambda^l(\Gamma^n(W)) \to \Gamma^n (\Lambda^l(W)),$$ $$ [w_1]_n \wedge \ldots \wedge [w_l]_n \mapsto [w_1 \wedge \ldots \wedge w_l]_n.$$Assume that $l < \dim(W)$  (which  does not at all follow from the proof above).  Then $\Psi$ is injective. Set $W':=\Lambda^l(W)$ and consider the composition of closed embeddings $$ \Gr(l,\Gamma^n(W)) \stackrel {\mathrm{Pl}}\hookrightarrow \P(\Lambda^l(\Gamma^n(W))) \stackrel {\Psi}\hookrightarrow \P(\Gamma^n (W')),$$ where $\mathrm{Pl}$ is the Pl\"ucker embedding. Denote by $L' \in \P(\Gamma^n(W'))$ the image of $L$ under this composite. One may then hope that, for a suitable choice of the data, the composite arrow  \[ G \stackrel \sim \to \Stab_{\PGL( W)}(L) \to \Stab_{\PGL( W')}(L')\] is an iso. 
\end{rem}

\begin{rem}\hfill \\
    A positive answer to Question \ref{Ql} would not simplify the proof of Theorem \ref{MainTh}.\\It may, however, be useful in other contexts.
\end{rem}
\section{Proof of Theorem \ref{MainTh}.}\label{SecProof}

Let $W$, $n$, $L$ and $l > \dim(W)$ be as in Proposition \ref{DividedTannaka}. Put $w:=\dim(W)$. Define \[V:=\Gamma^n(W)\] Denote by \[Z \subset  \P(V)\] the disjoint union of the closed subvarieties $\P(L)\simeq \P^{l-1}$  and $\Ver_n(\P(W))\simeq \P^{w-1}$. In the proof of Proposition \ref{DividedTannaka}, $W$ is fixed from the beginning, where $w$ can be picked arbitrarily large. The construction then works for all $n$ sufficiently large. It is straightforward to check that, when $n $ goes to infinity, so do $l$ and $\dim(V)-l$ (whereas $w$ stays fixed). In particular, one may assume that $w \neq l$ and $w \neq (\dim(V)-l)$.\\

 \begin{prop}\label{StabZ}
The natural inclusion $$G \hookrightarrow \Stab_{\PGL(V)}(Z \subset \P(V))$$  is an isomorphism of algebraic $F$-groups.
 \end{prop}
 
 \begin{dem}
 We first show that $$  \Stab_{\PGL(V)}(Z) =  \Stab_{\PGL(V)}(\Ver_n(\P(W))) \cap \Stab_{\PGL(V)}(\P(L)).$$ Inclusion $\supset$ is clear. To get equality, as both sides are linear algebraic groups over $F$, it suffices to prove  equality of their points, with values in a finite local  $F$-algebra $A$, which reads as\[  \Stab_{\PGL(V)}(Z)(A)= \Stab_{\PGL(V)}(\Ver_n(\P(W)))(A) \cap \Stab_{\PGL(V)}(\P(L))(A).\]
Pick $f \in \Stab_{\PGL(V)}(Z)(A)$. It induces an automorphism of the $A$-scheme $Z_A:=Z \times_F A$, which is the disjoint union of its irreducible clopen subschemes $\Ver(\P(W))_A$ and $\P(L)_A$. These are projective spaces of distinct dimensions, hence non-isomorphic. Thus $f$ preserves them both (which is a purely topological fact), proving the claim. To conclude,  apply Proposition  \ref{AutVeronese} combined to equality $G=\Stab_{\PGL(W)}(\P(L))$, provided by Proposition \ref{DividedTannaka}.
 \end{dem}
 
 Define  \[X:=\Bl_Z
 (\P(V)). \] 
The action of $G$ on $\P(V)$ stabilizes $Z$; hence an embedding of $F$-group schemes \[ \Phi: G \to \AAut(X) \]
 \begin{prop}
 The arrow $\Phi$ is an isomorphism.
 \end{prop}
 
 \begin{dem}
 We may assume $F=\overline F$.\\
 By Proposition \ref{StabZ}, we know that $G \stackrel \sim \to \Stab_{\AAut(\P(V))}(Z \subset \P(V))$.
 Using Proposition \ref{PropAutBl}, we thus know that $\Phi$ induces an iso \[ \mathbf J( G,\rho) \to  \mathbf J(\AAut(X),\rho), \] for every finite $F$-algebra $A$, with residue homomorphism $\rho: A \to F$.\\ To conclude, it remains to prove that  \[ \Phi(F): G(F) \to \AAut(X)(F) \] is onto, as a homomorphism of abstract groups. Denote by $E_1 \subset X$  (resp.  $E_2 \subset X$) the exceptional divisor lying over $ \P^{w-1}$  (resp. $ \P^{l-1}$). Since $w \neq l$ and $w \neq (\dim(V)-l)$,  Lemma \ref{LemEnd} implies that  $E_1$ and $E_2$ are non-isomorphic $F$-varieties.  Proposition \ref{PhiIso} then applies, concluding the proof.
 \end{dem}

\section{Acknowledgments.}
I am grateful to Michel Brion, for his reading and several decisive  exchanges.

\end{document}